\documentclass[12pt]{article}

\usepackage{amsmath,amsthm}
\usepackage{amssymb}
\usepackage{color}
\usepackage{xspace}
\usepackage[colorlinks=true,
linkcolor=green,
filecolor=brown,
citecolor=green]{hyperref}

\def\v{\vert}
\def\si{\sigma}
\def\a{\ensuremath{\mathcal A}\xspace}
\def\b{\ensuremath{\mathcal B}\xspace}
\def\c{\ensuremath{\mathcal C}\xspace}
\def\d{\ensuremath{\mathcal D}\xspace}
\def\e{\ensuremath{\mathcal E}\xspace}

\begin{document}
\newtheorem{lemma}{Lemma}
\newtheorem{theorem}{Theorem}
\newtheorem{prop}{Proposition}
\newtheorem*{cor}{Corollary}
\begin{center}
{\Large
Permutations Avoiding a Nonconsecutive Instance of a 2- or 3-Letter Pattern                           \\ 
}

\vspace{10mm}
DAVID CALLAN  \\
Department of Statistics  \\
\vspace*{-1mm}
University of Wisconsin-Madison  \\
\vspace*{-1mm}
1300 University Ave  \\
\vspace*{-1mm}
Madison, WI \ 53706-1532  \\
{\bf callan@stat.wisc.edu}  \\
\vspace{5mm}

\noindent November 1, 2006
\end{center}

\vspace{3mm}
\begin{abstract}
We count permutations avoiding a nonconsecutive instance of a two- or 
three-letter pattern, that is, the pattern may occur but only as consecutive 
entries in the permutation. Two-letter patterns give rise to the 
Fibonacci numbers. The counting sequences for the two representative 
three-letter patterns, 321 and 132, have respective generating functions
$(1+x^2)(C(x)-1)/(1+x+x^2-x C(x))$ and $C(x+x^{3})$ where $C(x)$ 
is the generating function for the Catalan numbers.
\end{abstract}

{\Large \textbf{1. Introduction}  } \quad
There is a large literature on pattern avoidance in permutations and 
words. Problems treated include counting permutations avoiding a 
pattern or set of patterns, or containing patterns a specified number 
of times. Or a pattern may be allowed but only as part of a 
larger pattern. Pattern occurrences may be unrestricted (the 
classical case),
or required to be consecutive (i.e., contiguous entries in the 
permutation), or a mixture of the two. Here we consider the following 
variation. Given a pattern $\pi$, how many permutations on $[n]$ avoid 
nonconsecutive instances of $\pi$, that is, the pattern $\pi$ may occur but only as consecutive 
entries in the permutation. In other words, \emph{all} occurrences of 
$\pi$ are required to be of the consecutive type. As usual, for 
$\pi=(\pi_{1},\pi_{2},\ldots,\pi_{k}$) a permutation on $[k]$, an instance of $\pi$ in a 
permutation $(a_{1},a_{2},\ldots,a_{n})$ is a subpermutation 
$(a_{i_{1}},a_{i_{2}},\ldots,a_{i_{k}}),\ 1\le i_{1}<\ldots<i_{k}\le 
n$, 
whose reduced form (replace 
smallest entry by 1, next smallest by 2, and so on) is $\pi$. 
In the present paper we deal with patterns of length $\le 3$.
The operations reverse and complement 
($a_{i}\to n+1-a_{i}$) on permutations and patterns preserve this type 
of pattern avoidance; so we need only 
consider representative patterns in the orbits under the action of 
the group they generate. For patterns of length $\le 3$, this reduces the problem 
(as usual) to the patterns 21, 321 and 132.
The following sections treat them one at a time.

\vspace{10mm}

{\Large \textbf{2. The 21 pattern}  }\quad
Avoiding a nonconsecutive 21 means that every inversion arises from a 
contiguous pair of entries.
Such permutations on $[n]$ are all obtained from the identity permutation 
$(1,2,\ldots,n) $ by selecting a set of non-overlapping pairs of 
consecutive integers in $[n]$ and, for each pair, switching the 
entries in those 2 positions. The first entries of these pairs form a 
scattered subset of $[n-1]$, that is, a subset for which distinct elements differ by at 
least 2. It is well known 
that such subsets are counted by the Fibonacci number $F_{n+1}$. 
Hence we have
\begin{theorem}
    The number of permutations on $[n]$ that avoid a nonconsecutive 
    $21$ pattern is $F_{n+1}$.
\end{theorem}

\vspace{10mm}

{\Large \textbf{3. The 321 pattern}  }\quad Let us introduce some 
notation:
\begin{itemize}
    \item  $\a_{n}$ is the full set of permutations on $[n]$ that 
    avoid a nonconsecutive 321.

    \item  $\b_{n},\d_{n}$ are, respectively, the permutations in 
    $\a_{n}$ that start with a 321 in the first 3 positions, and 
    those that don't. Thus $\v \a_{n}\v = \v \b_{n}\v + \v \d_{n}\v $.

    \item  $\c_{n}$ is the set of permutations in $\a_{n}$ that 
    contain no 321 at all, the set of so-called 321-avoiding 
    permutations. It is well known \cite{ec2} that $\v \c_{n}\v 
    =C_{n}$, the $n$th Catalan number.
    \item  $\a_{n,k},\ (1\le k \le n-2$) is the set of permutations 
    in $\a_{n}$ that contain one or more (necessarily consecutive) 321s, the first 
    of which starts at position $k$. Thus $\v \a_{n}\v  = 
    C_{n}+\sum_{k=1}^{n-2}\v \a_{n,k}\v$.   
\end{itemize}
First, observe that a permutation in $\b_{n}\ (n\ge 3)$ must have 2 as its 
second entry and 1 as its third entry (else an offending pattern would 
be present). Deleting these entries and 
subtracting 2 from all other entries is a bijection to $\d_{n-2}$. 
Hence $\v \b_{n}\v = \v \d_{n-2}\v$. Next, given $\rho \in \a_{n,k}\ 
(2\le k \le n-2)$, we can form two new permutations $\si,\tau$  as follows. 
Take the first $k-1$ entries and the entry in position $k+2$ 
(necessarily comprising the first $k$ positive integers) to form 
$\si$. Delete the first $k-1$ entries and reduce (replace 
smallest entry by 1, next smallest by 2, and so on) to form $\tau$. We 
leave the reader to verify that this is a bijection $\a_{n,k} \to 
\c_{k} \times \b_{n-k+1}$. Hence $\v \a_{n,k}\v = C_{k} \v 
\b_{n-k+1}\v = C_{k} \v \d_{n-k-1}\v$.

Now, with a lowercase letter $a_{n}$ denoting the size of $\a_{n}$ and so on, we have 
just shown that $a_{n}=C_{n}+\sum_{k=1}^{n-2}C_{k}d_{n-k-1}$ (since 
$C_{1}=1$). But also $a_{n}=b_{n}+d_{n}$ and $b_{n}=d_{n-2}$. Eliminating 
$a_{n}$ and $b_{n}$ and reindexing yields
\[
d_{n}=C_{n}+\sum_{k=1}^{n-3}C_{k+1}d_{n-2-k},
\]
a recurrence for $d_{n}$ that is valid for $n\ge 1$ and hence determines $d_{n}$ 
(no initial condition necessary).  
Since the sum is a convolution, this recurrence routinely 
yields the generating function
\[
D(x)=\frac{C^{*}(x)}{1+x^{2}-xC^{*}(x)},
\]
where $D(x):=\sum_{n\ge 1}d_{n}x^{n}$ and $C^{*}(x):=\sum_{n\ge 1}C_{n}x^{n}=
\frac{1-\sqrt{1-4x}}{2x}-1$ (the superscript $^{*}$ 
indicating that $C_{0}$ is omitted from the sum).
The formula for $a_{n}$ then yields 
\begin{theorem}
    The number $a_{n}$ of permutations on $[n]$ that avoid a nonconsecutive 
    $321$ pattern has generating function
    \[
\sum_{n\ge 1}a_{n}x^{n}=\frac{C^{*}(x)}{1-\frac{x}{1+x^{2}}C^{*}(x)} .
\]
\end{theorem}
The sequence $(a_{n})_{n\ge1}$ begins $1,\, 2,\, 6,\, 18,\, 56,\, 182,\, 607,\, 
2064,\ldots$ .
\vspace{10mm}

{\Large \textbf{3. The 132 pattern}  }\quad First, observe that 
132-avoiding permutations are characterized by the property that, for 
each entry $a$, the set of succeeding entries that are $<a$ form an 
initial segment of the positive integers. Also, Simion and Schmidt 
\cite{simion-schmidt} gave a bijection (see \cite{wilfequivalence} 
for an equivalent description) from 132-avoiding to 123-avoiding 
permutations, the latter corresponding to 321-avoiding under reversal, 
and so 132-avoiding permutations are also counted by the Catalan 
numbers.

Now let $\e_{n}$ denote the set of permutations on $[n]$ that 
avoid a nonconsecutive 132 and $\e_{n,k}$ the permutations in 
$\e_{n}$ with $k$ (necessarily consecutive) 132s. Consider a permutation 
in $\e_{n,k}$. First, record the positions $i_{1},i_{2},\ldots,i_{k}$ 
of the middle entries of its 132 patterns---they form a subset of the interval $[2,n-1]$ 
whose entries all differ by at least 3 because, as is easily checked, 
the 132 patterns cannot overlap. Then delete the first and last entry 
of each 132 pattern and reduce. The result is a 132-avoiding permutation 
on $[n-2k]$. Thus each permutation in $\e_{n,k}$ yields a set $\{ 
i_{1},i_{2},\ldots,i_{k}\} \subset [2,n-1]$ with entries differing by 
at least 3 together with a 132-avoiding permutation on $[n-2k]$.

Conversely, given any such set and permutation, we can construct a 
permutation in $\e_{n,k}$ that yields them.
Take, for example, $n=10$ and $k=2$ and suppose given $\{ 
i_{1},i_{2},\ldots,i_{k}\}=\{4,8\}$ and 132-avoiding permutation 
$(6,5,3,4,2,1)$. Place blanks in the positions immediately 
neighboring each $i_{j}$ and place the permutation entries in the 
remaining positions:
\[
    \begin{array}{cccccccccc}
        1 & 2 & 3 & 4 & 5 & 6 & 7 & 8 & 9 & 10  \\
        6 & 5 &  & 3 &  & 4 &  & 2 &  & 1
    \end{array}
\]
Leave intact the entries after the entry in position $i_{1}$ and 
smaller than it (as noted above, these entries form an initial 
segment $1,2,\ldots,j$
of the positive integers), place $j+1$ in position $i_{1}-1$, $j+2$ 
in position $i_{1}+1$, and increase all other entries by $2$:
\[
    \begin{array}{cccccccccc}
        1 & 2 & 3 & 4 & 5 & 6 & 7 & 8 & 9 & 10  \\
        8 & 7 & 3 & 5 & 4 & 6 &  & 2 &  & 1
    \end{array}
\]
Proceed likewise left to right to fill in the remaining pairs of  blanks. 
It is easy to check that the result, here $(10,9,5,7,6,8,2,4,3,1)$,
is a permutation in $\e_{n,k}$ that yields the given $\{ 
i_{1},i_{2},\ldots,i_{k}\}$ and 132-avoiding permutation under the 
process of the previous paragraph.

There are $\binom{n-2k}{k}$ such sets $\{ 
i_{1},i_{2},\ldots,i_{k}\}$ and $C_{n-2k}$ 132-avoiding permutations 
on $[n-2k]$ and so $\v \e_{n,k} \v=\binom{n-2k}{k}C_{n-2k}$. 
Summing over $k$ and taking into account the permutations that avoid 132 
altogether, we have
\begin{theorem}
The number of permutations on $[n]$ that avoid a nonconsecutive 
$132$ pattern is
\[
\sum_{k=0}^{\lfloor n/3 
\rfloor}\binom{n-2k}{k}C_{n-2k}.
\]
\end{theorem}
Routine manipulations lead to a succinct generating function:
\begin{eqnarray*}
    \sum_{n\ge 0}\v \e_{n}\v x^{n}& = & \sum_{n,k \ge 0} 
    \binom{n-2k}{k}C_{n-2k}x^{n} \\
     & = & \sum_{k \ge 0}\bigg(\sum_{n\ge 3k}\binom{n-2k}{k}C_{n-2k}x^{n-3k}\bigg)x^{3k}  \\
      & = & \sum_{k \ge 0}\bigg(\sum_{n\ge k}\binom{n}{k}C_{n}x^{n-k}\bigg)x^{3k}  \\
     & = &  \sum_{k \ge 0} \frac{C^{\{k\}}(x)}{k!}x^{3k} \\
     & = & C(x+x^{3}),
\end{eqnarray*}
the last equality by Taylor's theorem, where $C(x)=\frac{1-\sqrt{1-4x}}{2x}$ is the 
generating function for the Catalan numbers. So we have
\begin{cor}
    The generating function for permutations that avoid a nonconsecutive 
$132$ pattern is 
\[
\frac{1-\sqrt{1-4x-4x^{3}}}{2(x+x^{3})}.
\]
\end{cor}
The counting sequence ($n\ge 0$) begins  $1,\,1,\, 2,\, 6,\, 18,\, 57,\, 190,\, 654,\, 
2306,\ldots$ .

\textbf{Added in Proof}\quad These results have previously been 
obtained in a wider context by Anders Claesson \cite{segmented}.


\begin{thebibliography}{99}
\bibitem{ec2} Richard P.~Stanley, \emph{Enumerative Combinatorics} 
Vol.\,2, Cambridge University Press, 1999. Exercise 6.19 and related 
material on Catalan numbers are available 
online at
\htmladdnormallink{http://www-math.mit.edu/$\,\widetilde{\ }\,$rstan/ec/ }{http://www-math.mit.edu/~rstan/ec/}.

\bibitem{simion-schmidt}
Rodica Simion and Frank W. Schmidt, 
Restricted Permutations, 
\emph{Europ. J. Combinatorics} {\bf 6} (1985), 383-406. 

\bibitem{wilfequivalence} David Callan, A Wilf Equivalence Related to Two Stack 
Sortable Permutations, 2005, preprint,
\htmladdnormallink{http://front.math.ucdavis.edu/math.CO/0510211 }{http://front.math.ucdavis.edu/math.CO/0510211 }.

\bibitem{segmented}Anders Claesson, Counting segmented permutations using bicoloured Dyck paths, 
\htmladdnormallink{\emph{Elec. J. 
Comb.}}{http://www.combinatorics.org/} ,Vol. 12, 
\textbf{R39}, 2005.

\end{thebibliography}
\end{document}